\documentclass[12pt,a4paper]{article}
\usepackage{theorem}
\usepackage{amsmath,amssymb,amscd,latexsym}

\newtheorem{thm}{Theorem}[section]

\newtheorem{prop}[thm]{Proposition}
\newtheorem{cor}[thm]{Corollary}
\newtheorem{defn}[thm]{Definition}
{\theorembodyfont{\rmfamily}  }

\newcommand{\Dom}{\operatorname{Dom}}
\newcommand{\axis}{\operatorname{axis}}

\newcommand{\Stab}{\operatorname{Stab}}

\begin{document}

\title{\Large{\textbf{Strongly singular MASA's and mixing actions in
      finite von Neumann algebras}}}
\author{Paul Jolissaint and Yves Stalder}
\maketitle
\begin{abstract}
Let $\Gamma$ be a countable group and let
$\Gamma_0$ be an infinite abelian subgroup of $\Gamma$. We prove
that if the pair $(\Gamma,\Gamma_0)$ satisfies 
some combinatorial condition called (SS), then the abelian
subalgebra $A=L(\Gamma_0)$ is a singular MASA in
$M=L(\Gamma)$ which satisfies a weakly mixing condition. 
If moreover it satisfies a stronger condition called (ST), then it
provides a singular MASA with a strictly stronger mixing property.
We describe families of examples of both types coming from free
products, HNN extentions and semidirect products, and in particular
we exhibit examples of singular MASA's that
satisfy the weak mixing condition but not the
 strong mixing one.
\par\vspace{3mm}\noindent
\end{abstract}
\section{Introduction}

If $M$ is a von Neumann algebra, if $A$ is a 
maximal abelian von
Neumann subalgebra of $M$ (henceforth called a MASA), let $\mathcal{N}_M(A)$ be the
\emph{normaliser} of $A$ in $M$: it is the subgroup of the unitary
group $U(M)$ of all elements $u$ such that $uAu^*=A$. Then $A$ is
\emph{singular} in $M$ if $\mathcal{N}_M(A)$ is as small as possible,
namely, if $\mathcal{N}_M(A)=U(A)$.
Until recently, it 
was quite difficult in general to exhibit
singular MASA's in von Neumann algebras, though S. Popa proved among others in
\cite{Popa}
that all separable type $\textrm{II}_1$ factors admit singular MASA's.
\par\vspace{3mm}\noindent
\textbf{Example.} This example is due to F. Radulescu \cite{Rad}. Let $L(F_N)$ be the factor
associated to the non abelian
free group on $N$ generators $X_1, \ldots, X_N$ and let $A$ be the
abelian von Neumann subalgebra generated by $X_1+ \ldots
+X_N+X_1^{-1}
+ \ldots  +X_N^{-1}$. Then $A$ is a singular MASA in $L(F_N)$.
$A$ is called the \emph{radial} or \emph{Laplacian} subalgebra because
its elements coincide with convolution operators by functions that
depend only on the length of the elements of $F_N$.
\par\vspace{3mm}
Recently, T. Bildea generalized Radulescu's example in
\cite{Bil} by using the notion of asymptotic homomorphism (see below): 
Let $G$ be
either $F_N$ or a free product $G_1\star\ldots\star G_m$ of $m\geq 3$
groups, each finite of order $p\geq 2$ with $m\geq p$. Then the radial
subalgebra is a singular MASA in $L(G)^{\bar{\otimes}_k}$. 
\par\vspace{3mm}
Motivated by S. Popa's articles \cite{Popa} and \cite{Popa83}, the authors of 
 \cite{RSS} and \cite{SS}
introduced  sufficient conditions on an abelian von
Neumann subalgebra $A$ of a finite von Neumann algebra $M$ that
imply that $A$ is even a strongly
singular MASA in $M$. This means that $A$ satisfies the apparently stronger
condition: 
for all $u\in U(M)$, one has
$$
\sup_{x\in M,\Vert x\Vert\leq 1}\Vert E_A(x)-E_{uAu^*}(x)\Vert_2\geq
\Vert u-E_A(u)\Vert_2.
$$
In fact, it was proved by A. Sinclair, R. Smith, S. White and A. Wiggins in
\cite{SSWW}
that all singular MASA's are strongly singular. Nevertheless,
it is sometimes easier to prove directly strong singularity.
\begin{prop}(\cite{RSS})
Suppose that the pair $A\subset M$ satisfies the following condition 

(WM):
$\forall x,y\in M$ and $\forall \varepsilon>0$, there exists $v\in
U(A)$ 
such that
$$
\Vert E_A(xvy)-E_A(x)vE_A(y)\Vert_2\leq \varepsilon.
$$
Then $A$ is a strongly singular MASA in $M$.
\end{prop}
Let us reproduce an adaptation of the proof of Lemma 2.1 of \cite{RSS}
for convenience:
\par\vspace{3mm}\noindent
\emph{Proof.} Fix $u\in U(M)$ and $\varepsilon>0$, and take $x=u^*$,
$y=u$. There exists $v\in U(A)$ such that 
$$
\Vert E_A(u^*vu)-E_A(u^*)vE_A(u)\Vert_2=
\Vert E_A(v^*u^*vu)-E_A(u^*)E_A(u)\Vert_2\leq \varepsilon.
$$
(Commutativity of $A$ is crucial here.)
Hence, we get :
\begin{eqnarray*}
  \Vert E_A-E_{uAu^*}\Vert_{\infty,2}^2 & \geq &
\Vert v-uE_A(u^*vu)u^*\Vert_2^2\\
& = & \Vert u^*vu-E_A(u^*vu)\Vert_2^2\\
& = & 1-\Vert E_A(u^*vu)\Vert_2^2\\
& \geq & 1-(\Vert E_A(u^*)vE_A(u)\Vert_2+\varepsilon)^2\\
& \geq & 1- (\Vert E_A(u)\Vert_2+\varepsilon)^2\\
& = & \Vert u-E_A(u)\Vert_2^2-2\varepsilon\Vert E_A(u)\Vert_2-\varepsilon^2.
\end{eqnarray*}
As $\varepsilon$ is arbitrary, we get the desired inequality.
\hfill $\square$
\par\vspace{3mm}
Earlier, in \cite{SS}, A. Sinclair and R. Smith used a stronger condition
(that we call \emph{condition (AH)} here) in order to get singular MASA's:
\par
Given $v\in U(A)$, the conditional expectation $E_A$ is an
\emph{asymptotic homomorphism with respect to} $v$ if 
$$
(AH)\quad\lim_{|k|\to\infty}\Vert E_A(xv^ky)-E_A(x)v^kE_A(y)\Vert_2=0
$$
for all $x,y\in M$.
\par\vspace{3mm}
Both conditions (WM) and (AH) remind one of mixing properties of group
actions, because of the following equality (since $A$ abelian):
$$
\Vert E_A(vxv^*y)-E_A(x)E_A(y)\Vert_2=\Vert
E_A(xv^*y)-E_A(x)v^*E_A(y)\Vert_2
$$ 
$\forall x,y\in M,\ \forall v\in U(A).$ 
\par\vspace{3mm}
Thus our point of view here is the following: every subgroup $G$ of $U(A)$ 
acts by inner automorphisms on $M$ by conjugation
$$
\sigma_v(x)=vxv^*\quad\forall v\in G,\ \forall x\in M,
$$
and as will be seen, condition (WM) is in some sense
equivalent to a weakly mixing action of $G$, and condition (AH) is
equivalent to a strongly mixing action of the cyclic subgroup of
$U(A)$ generated by the distinguished unitary $v$.
\par\vspace{3mm}
Section 2 is devoted to weakly mixing actions of subgroups $G\subset
U(A)$, to crossed products by weakly mixing actions of countable groups,
and to pairs $(\Gamma,\Gamma_0)$ where $\Gamma_0$ is an abelian
subgroup of the countable group $\Gamma$ which provide pairs
$A=L(\Gamma_0)\subset M=L(\Gamma)$ with $A$ weakly mixing in $M$ (see
Definition 2.1).
It turns out that the weak mixing property is completely determined by
a combinatorial property of the pair $(\Gamma,\Gamma_0)$ 
called \emph{condition (SS)},
already introduced in \cite{RSS} to provide examples of strongly
singular MASA's. 
\par\vspace{3mm} 
Section 3 is devoted to a (strictly) stronger condition 
(called \emph{condition (ST)})
on pairs 
$(\Gamma,\Gamma_0)$ as above  
which is related to strongly mixing
actions of groups. 
\par\vspace{3mm}
Section 4 contains various families of examples (free products with
amalgamation, HNN extentions, semidirect products), and some of them
 prove
that condition (ST) is strictly stronger than condition (SS) so that 
they provide 
two distinct levels of ``mixing MASA's'': the weak ones and the
strong ones.
\par\vspace{3mm}
We will see that
many pairs $(\Gamma,\Gamma_0)$ satisfy a strictly stronger condition than 
condition (ST), based on malnormal subgroups: 
 $\Gamma_0$ is said to be a 
\emph{malnormal subgroup} of $\Gamma$ if
for every $g\in\Gamma\setminus\Gamma_0$,
  one has $g\Gamma_0g^{-1}\cap\Gamma_0=\{1\}$.
Such pairs have been considered first in the pioneering article
\cite{Popa83} in order
to control normalizers (in particular relative commutants) of
$L(\Gamma_0)$ and of its diffuse subalgebras,
and, as a byproduct, to produce singular MASA's $L(\Gamma_0)$.
They were also used in 
\cite{Rob},
\cite{SS} to provide more examples of (strongly) singular MASA's in type
$\textrm{II}_1$ factors that fit Popa's criteria of Proposition 4.1 in
\cite{Popa83}. 

\par\vspace{3mm}
\emph{Acknowledgements.} We are grateful to A. Sinclair and A. Valette
for helpful
comments, and to S. Popa
for having pointed out relationships between our article and relative
mixing conditions appearing in 
\cite{PopaOA} and \cite{PopaGR}, and the use of malnormal subgroups in
\cite{Popa83}.


\section{Weak mixing}

In the rest of the article, $M$ denotes a finite von Neumann algebra,
and $\tau$ is some normal, faithful, finite, normalised trace on $M$
(henceforth simply called a \emph{trace} on $M$).
It defines a scalar product on $M$: $\langle
a,b\rangle=\tau(b^*a)=\tau(ab^*)$, and the corresponding completion is
the Hilbert space $L^2(M,\tau)$ on which $M$ acts by left
multiplication extending the analogous operation on $M$. 
As usual, we denote by $\Vert\cdot\Vert_2$
the corresponding Hilbert norm. When the trace $\tau$ must be
specified, we write $\Vert\cdot\Vert_{2,\tau}$.
We denote
also by $M_\star$ the predual of $M$, i.e. the Banach space of all
normal linear functionals on $M$. We will always assume for
convenience that $M_{\star}$ is separable, or equivalently, that
$L^2(M,\tau)$ is a separable Hilbert space.
Recall that, for every $\varphi\in
M_\star$, there exist $\xi,\eta\in L^2(M,\tau)$ such that
$\varphi(x)=\langle x\xi,\eta\rangle$ for all $x\in M$.
\par\vspace{3mm}
Let $\Gamma$ be a countable group and let $\Gamma_0$ be an abelian
subgroup of $\Gamma$. Denote by $L(\Gamma)$ (respectively $L(\Gamma_0)$) the von
Neumann algebra 
generated by the left regular representation $\lambda$
of $\Gamma$ (respectively $\Gamma_0$). Recall that
$\lambda:\Gamma\rightarrow U(\ell^2(\Gamma))$ is defined by
$(\lambda(g)\xi)(h)=\xi(g^{-1}h)$ for all $g,h\in \Gamma$ and
$\xi\in\ell^2(\Gamma)$. It extends linearly to the group algebra
$\mathbb{C}\Gamma$, and $L(\Gamma)$ is the weak-operator closure of
$\lambda(\mathbb{C}\Gamma)=:L_f(\Gamma)$. 
The normal functional
$\tau(x)=\langle x\delta_1,\delta_1\rangle$ is a faithful trace on $L(\Gamma)$.
For $x\in L(\Gamma)$, denote by
$\sum_{g\in\Gamma}x(g)\lambda(g)$ its ``Fourier expansion'':
$x(g)=\tau(x\lambda(g^{-1}))$ for every $g\in \Gamma$, and
the series $\sum_{g}x(g)\lambda(g)$ converges to $x$ in the
$\Vert\cdot\Vert_2$-sense so that
 $\sum_{g\in
  \Gamma}\vert x(g)\vert^2=\Vert x\Vert_2^2$. 
\par\vspace{3mm}
Let $1\in B$ be a unital von Neumann subalgebra of 
the von Neumann algebra
 $M$ gifted with some trace $\tau$ as above and let $E_B$
be the $\tau$-preserving conditional expectation of $M$ onto $B$. Then
$E_B$ is characterised by
the following two conditions: $E_B(x)\in B$ for all $x\in M$ and
$\tau(E_B(x)b)=\tau(xb)$ for all $x\in M$ and all $b\in B$.
It enjoys the well-known properties:
\begin{enumerate}
\item [(1)] $E_B(b_1xb_2)=b_1E_B(x)b_2$ for all $x\in M$ and all
  $b_1,b_2\in B$;
\item [(2)] $\tau\circ E_B=\tau$.
\item [(3)] If $M=L(\Gamma)$ is the von Neumann algebra associated to the
  countable group $\Gamma$, if
$H$ is a subgroup of $\Gamma$, if $B=L(H)$ and if 
$x\in M$, then  $E_B(x)=\sum_{h\in H}x(h)\lambda(h)$.
\end{enumerate} 
Let also $\alpha$ be a
$\tau$-preserving action of $\Gamma$ on $M$. Recall that it is
\emph{weakly mixing} if, for every finite set $F\subset M$ and for
every $\varepsilon>0$, there exists $g\in\Gamma$ such that
$$
|\tau(\alpha_g(a)b)-\tau(a)\tau(b)|<\varepsilon\quad\forall a,b\in F.
$$
\par\vspace{3mm}
In \cite{PopaGR}, S. Popa introduced a relative version of weakly
mixing actions:\\
If $1\in A\subset M$ is a von Neumann subalgebra such that $\alpha_g(A)=A\
\forall g\in\Gamma$, the action $\alpha$ is called \emph{weakly mixing
  relative to} $A$ if, for every finite set $F\subset M\ominus
A:=\{x\in M\ :\ E_A(x)=0\}$, for
every $\varepsilon>0$, one can find $g\in\Gamma$ such that 
$$
\Vert E_A(x^*\alpha_g(y))\Vert_2\leq \varepsilon\quad\forall x,y\in F.
$$
Here we consider a (countable) subgroup $G$ of the
unitary group $U(A)$ which acts on $M$ by conjugation:
$$
\sigma_v(x)=vxv^*\quad\forall v\in G,\ \forall x\in M.
$$
This allows us to introduce the following definition:

\begin{defn}
The abelian von Neumann subalgebra 
$A$ is \textbf{weakly mixing in} $M$ if there
exists a subgroup $G$ of $U(A)$ such that the corresponding action by
conjugation is weakly mixing relative to $A$ in Popa's sense.  
\end{defn}
Notice that, $A$ being abelian, it is equivalent to asking that, for every finite set
$F\subset M$ and every $\varepsilon>0$, there exists $v\in G$ such
that 
$$
\Vert E_A(xvy)-E_A(x)vE_A(y)\Vert_2\leq \varepsilon\quad\forall x,y\in
F.
$$
Weakly mixing actions provide singular MASA's, as was already kwown to
many people; see for instance \cite{NS}, \cite{PopaOA} and
\cite{RSS}. 
More precisely, one has:

\begin{prop}
Let $\Gamma_0$ be an abelian group which acts
on a finite von Neumann algebra $B$ and which preserves a trace
$\tau$, then the abelian von Neumann subalgebra $A=L(\Gamma_0)$ of the
crossed product $M=B\rtimes\Gamma_0$ is weakly mixing in $M$ if and
only if the
action of $\Gamma_0$ is.
\end{prop}
\emph{Sketch of proof.} Observe first that, for $x=\sum_\gamma
x(\gamma)\lambda(\gamma)\in M$, the projection of $x$ onto $A$ is
given by $E_A(x)=\sum_\gamma \tau(x(\gamma))\lambda(\gamma)$. 
\par
Suppose that the action $\alpha$ of $\Gamma_0$ is weakly
mixing. It suffices to prove that, for all finite sets
$E\subset\Gamma_0$ and $F\subset B$ and for every $\varepsilon>0$,
there exists $\gamma\in\Gamma_0$ such that
$$
\Vert
E_A(a\lambda(g)\lambda(\gamma)b\lambda(h))-E_A(a\lambda(g))\lambda(\gamma)E_A(b\lambda(h))
\Vert_2<\varepsilon
$$
for all $a,b\in F$ and all $g,h\in E$. But,
$E_A(a\lambda(g))=\tau(a)\lambda(g)$,
for all $a\in F$ and $g\in \Gamma$, and 
$E_A(a\lambda(g)\lambda(\gamma)b\lambda(h))=\tau(a\alpha_\gamma(\alpha_g(b)))\lambda(g\gamma
h)$
which suffices to get the conclusion.
\par
Conversely, if $G$ is a subgroup of $U(A)$ whose action is weakly
mixing relative to $A$, let $F$ be as above and let
$\varepsilon>0$. There exists a unitary $u\in G$ such that 
$$
\sum_{a,b\in F}\Vert
E_A(aub)-E_A(a)uE_A(b)\Vert_2^2<\frac{\varepsilon^2}{2}.
$$
As above, since $u=\sum_\gamma u(\gamma)\lambda(\gamma)$ with
$u(\gamma)\in \mathbb C$ for every $\gamma$, one has
$$
E_A(aub)=\sum_\gamma u(\gamma)\tau(a\alpha_\gamma(b))\lambda(\gamma)
$$
and 
$$
E_A(a)uE_A(b)=\sum_\gamma
u(\gamma)\tau(a)\tau(b)\lambda(\gamma).
$$
This implies that 
$$
\sum_{a,b\in F}\sum_\gamma |u(\gamma)|^2|\tau(a\alpha_\gamma(b))-\tau(a)\tau(b)|^2<
\frac{\varepsilon^2}{2}.
$$
As $\sum_\gamma |u(\gamma)|^2=1$, this implies the existence of
$\gamma\in\Gamma_0$ such that
$$
\sum_{a,b\in
  F}|\tau(a\alpha_\gamma(b))-\tau(a)\tau(b)|^2<\varepsilon^2.
$$
\hfill $\square$
\par\vspace{3mm}
When $M=L(\Gamma)$ and $A=L(\Gamma_0)$, where $\Gamma_0$ is an abelian
subgroup of $\Gamma$, it turns out that weak mixing of $A$ is
equivalent to a combinatorial property of the pair of groups
$(\Gamma,\Gamma_0)$ as it appears in \cite{RSS}:

\begin{prop}
For a pair $(\Gamma,\Gamma_0)$ as above, the following two conditions
are equivalent:
\begin{enumerate}
\item [(1)] (SS) For every finite subset $C\subset
  \Gamma\setminus\Gamma_0$, there exists $\gamma\in\Gamma_0$ such that
  $g\gamma h\notin\Gamma_0$ for all $g,h\in C$;
\item  [(2)] (WM) $A=L(\Gamma_0)$ is weakly mixing in $M=L(\Gamma)$. 
\end{enumerate}
\end{prop}
\emph{Proof.} If $(\Gamma,\Gamma_0)$ satisfies condition (SS), take
$G=\lambda(\Gamma_0)$. Let $F$ be a finite subset of $L(\Gamma)$ such
that $E_A(x)=0$ for every $x\in F$. By density, one assumes that $F$
is contained in $L_f(\Gamma)$. We intend to prove that
there exists $\gamma\in\Gamma_0$ such that $E_A(x\lambda(\gamma)y)=0$
for all $x,y\in F$. By the above assumptions, there exists a finite
set $C\subset \Gamma\setminus\Gamma_0$ such that every $x\in F$ has
support in $C$. Fix $x$ and $y$ in $F$ and choose $\gamma\in\Gamma_0$
as in condition (SS) with respect to $C$. Then
\begin{eqnarray*}
  E_A(x\lambda(\gamma)y) & = & E_A\left(
\sum_{g,h\in C}x(g)y(h)\lambda(g\gamma h)\right)\\
 & = & \sum_{g,h\in C,g\gamma h\in\Gamma_0} x(g)y(h)\lambda(g\gamma h)=0
\end{eqnarray*}
since the set $\{(g,h)\in C\times C\ :\ g\gamma h\in\Gamma_0\}$ is
empty. This shows that $A$ is weakly mixing in $M$.
\par
Conversely, if $A$ satisfies the latter condition in $M$, choose a subgroup $G$ of
$U(A)$ such that its action on $M$ is weakly mixing relative to $A$
and let $C\subset \Gamma\setminus\Gamma_0$ be a finite set. Take
$F=\lambda(C)$, and observe that $E_A(\lambda(g))=0$ for every $g\in
C$. Finally, choose any $0<\varepsilon<1/2$. Then there exists $v\in
G$ such that 
$$
\sum_{g,h\in C}\Vert E_A(\lambda(g)v\lambda(h)\Vert^2_2<\varepsilon^2.
$$
We have $v=\sum_{\gamma\in\Gamma_0}v(\gamma)\lambda(\gamma)$ and 
$\Vert v\Vert_2^2=1=\sum_\gamma |v(\gamma)|^2$. Moreover,
\begin{eqnarray*}
  E_A(\lambda(g)v\lambda(h)) & = & 
\sum_{\gamma\in\Gamma_0}v(\gamma)E_A(\lambda(g\gamma h))\\
& = & \sum_{\gamma\in\Gamma_0,g\gamma
  h\in\Gamma_0}v(\gamma)\lambda(g\gamma h).
\end{eqnarray*}
Thus,
$$
\Vert E_A(\lambda(g)v\lambda(h)\Vert^2_2=
\sum_{\gamma\in\Gamma_0,g\gamma
  h\in\Gamma_0}|v(\gamma)|^2.
$$
Assume that for every $\gamma\in\Gamma_0$ one can find $g_\gamma$ and
$h_\gamma$ in $C$ such that $g_\gamma\gamma h_\gamma\in
\Gamma_0$. Then we would get 
\begin{eqnarray*}
  1=\sum_{\gamma\in\Gamma_0}|v(\gamma)|^2 & =& 
\sum_{\gamma\in\Gamma_0,g_\gamma\gamma
  h_\gamma\in\Gamma_0}|v(\gamma)|^2\\
& \leq & \sum_{g,h\in C}\sum_{\gamma\in\Gamma_0,g\gamma h\in\Gamma_0}|v(\gamma)|^2<\varepsilon^2
\end{eqnarray*}
which is a contradiction. 
\hfill $\square$
\par\vspace{3mm}
Before ending the present section, let us recall  
examples of pairs $(\Gamma,\Gamma_0)$ that satisfy
condition (SS). They are
taken from \cite{RSS}. Let $\Gamma$ be a group of isometries of some
metric space $(X,d)$ and let $\Gamma_0$ be an abelian subgroup of
$\Gamma$.
Assume that there is a
 $\Gamma_0$-invariant subset
$Y$ of $X$ such that
\begin{enumerate}
\item [(C1)] there exists a compact set $K\subset Y$ such that
  $\Gamma_0 K=Y$ ;
\item [(C2)] if $Y\subset_\delta g_1Y\cup g_2Y\cup\ \ldots\ \cup g_nY$ for some
  $g_j$'s in $\Gamma$, and some $\delta>0$, then there exists $j$ such
  that $g_j\in\Gamma_0$. (Recall that for subsets $P,Q$ of $X$ and $\delta>0$, $P\subset_\delta
Q$ means that $d(p,Q)\leq\delta$ for every $p\in P$.)
\end{enumerate}

Then it is proved in Proposition 4.2 of \cite{RSS} that 
the pair $(\Gamma,\Gamma_0)$ satisfies condition (SS).
\par\vspace{3mm}
Now let $G$ be a semisimple Lie group
with no centre and no compact factors. Let $\Gamma$ be a torsion free
cocompact lattice in $G$. Then $\Gamma$ acts freely on the symmetric space
$X=G/K$, where $K$ is a maximal compact subgroup of $G$, and the
quotient manifold $M=\Gamma\backslash X$ has universal covering space
$X$.
In particular, $\Gamma$ is the fundamental group $\pi(M)$ of $M$. Let
$r$ be the rank of $X$ and let
 $T^r\subset M$ be a totally geodesic embedding of a flat
$r$-torus in $M$, so that the inclusion $i:T^r\rightarrow M$ induces
an injective homomorphism $i_*:\pi(T^r)\rightarrow \pi(M)$. Thus
$\Gamma_0=i_*\pi(T^r)\cong \mathbb{Z}^r$ is an abelian subgroup of
$\Gamma$ in a natural way. Then it is proved in Theorem 4.9 of
\cite{RSS} that the pair $(\Gamma,\Gamma_0)$ satisfies conditions (C1)
and (C2) above, hence that it satisfies condition (SS) as well. In the
same vein, the authors also get examples of pairs coming from groups
acting cocompactly on locally finite euclidean buildings.


\section{Strong mixing}

Let $F$ be the Thompson's group; it admits the following
presentation:
$$
F=\langle x_0,\ x_1,\ \ldots\ |\ x_i^{-1}x_nx_i=x_{n+1},\ 0\leq
i<n\rangle.
$$ 
Let $\Gamma_0$ be the subgroup generated by $x_0$.
In \cite{JOF}, the first named author proved that
the pair $(F,\Gamma_0)$
 satisfies a stronger property than condition (SS) that was called \emph{condition
  (ST)} and which is described as follows: 

\begin{defn}
Let $(\Gamma,\Gamma_0)$ be a pair as above. 
Then we say that it satisfies \textbf{condition (ST)} if,
for every finite subset $C\subset\Gamma\setminus\Gamma_0$
  there exists a finite subset $E\subset\Gamma_0$ such that
  $gg_0h\notin\Gamma_0$ for all $g_0\in\Gamma_0\setminus E$ and all
$g,h\in C$.
\end{defn}
\par\vspace{3mm}\noindent
\textbf{Remark.}  Observe that, taking finite unions of exceptional
sets $E$ of $\Gamma_0$ if necessary, condition (ST) is equivalent to:
\par
\emph{For all $g,h\in\Gamma\setminus
  \Gamma_0$, there exists a finite subset $E$ of $\Gamma_0$ such that
$g\gamma h\notin\Gamma_0$ for all
  $\gamma\in \Gamma_0\setminus E$.}
\par
Since the subset $\Gamma\setminus\Gamma_0$ is stable under the mapping
$g\mapsto g^{-1}$,
when $\Gamma_0$ is an infinite cyclic group generated by some element
$t$, condition (ST) is still equivalent to:
\par
\emph{For all $g,h\in \Gamma\setminus \Gamma_0$, there exists a
  positive integer $N$ such that, for every $|k|> N$, one has
$gt^kh\notin\Gamma_0$.}
\par\vspace{3mm}
For future use in the present section, for every subset $S$ of a group
$\Gamma$ we put $S^*=S\setminus\{1\}$.
\par\vspace{3mm}
We observe in the next section that condition (ST) is strictly stronger
than condition (SS); examples are borrowed from Section 5 of
\cite{SSP}.
As it is the case for condition (SS), 
it turns out that condition (ST) is completely characterized by the
pair of von Neumann algebras $L(\Gamma_0)\subset L(\Gamma)$.
\par\vspace{3mm}
 To prove that, we need some facts and definitions.
Let $M$ and $\tau$ be as above. Let us say that a subset $S$
 of $M$ is \emph{$\tau$-orthonormal} if
$\tau(xy^*)=\delta_{x,y}$ for all $x,y\in S$.
We will need a weaker notion which is independent of the chosen trace $\tau$.

\begin{prop}
Let $M$ be a finite von Neumann algebra, let $\tau$ be a finite trace
on $M$ as above and let $S$ be an infinite subset of the unitary group
$U(M)$. The following conditions are equivalent:
\begin{enumerate}
\item [(1)] for every $\varphi\in M_\star$ and for every $\varepsilon>0$,
  there exists a finite subset $F$ of $S$ such that
  $|\varphi(u)|\leq\varepsilon$ for all $u\in S\setminus F$;
\item [(2)] for every $x\in M$ and for every $\varepsilon>0$, there
  exists a finite set $F\subset S$ such that $|\tau(ux)|\leq
  \varepsilon$ for all $u\in S\setminus F$;
\item [(2')] for any trace $\tau'$ on $M$, 
for every $x\in M$ and for every $\varepsilon>0$, there
  exists a finite set $F\subset S$ such that $|\tau'(ux)|\leq
  \varepsilon$ for all $u\in S\setminus F$;
\item [(3)] for every $\tau$-orthonormal finite set
  $\{x_1,\ldots,x_N\}\subset M$ and for every $\varepsilon>0$ there
  exists a finite set $F\subset S$ such that
$$
\sup\{|\tau(ux^*)|\ ;\ 
x\in \mathrm{span}\{x_1,\ldots,x_N\},\ \Vert x\Vert_2\leq
1\}\leq\varepsilon\quad\forall u\in S\setminus F;
$$
\item [(3')] for every trace $\tau'$ on $M$, for every $\tau'$-orthonormal finite set
  $\{x_1,\ldots,x_N\}\subset M$ and for every $\varepsilon>0$ there
  exists a finite set $F\subset S$ such that
$$
\sup\{|\tau'(ux^*)|\ ;\ 
x\in \mathrm{span}\{x_1,\ldots,x_N\},\ \Vert x\Vert_{2,\tau'}\leq
1\}\leq\varepsilon\quad\forall u\in S\setminus F;
$$ 
\end{enumerate}
In particular, if $S\subset U(M)$ satisfies the above conditions, if
$\theta$ is a $*$-isomorphism of $M$ onto some von Neumann algebra
$N$, then $\theta(S)\subset U(N)$ satisfies the same conditions.
\end{prop}
\emph{Proof.} $(1)\ \Rightarrow\ (2')\ \Rightarrow\ (2)$  
and $(3')\ \Rightarrow\ (3)$
are trivial.\\
$(2)\ \Rightarrow\ (3')$: If $\tau'$ is a trace on $M$, if
$\{x_1,\ldots,x_N\}\subset M$ is
$\tau'$-orthonormal and if $\varepsilon>0$ is fixed, there exists 
$h\in M$ such that 
$$
\Vert \tau'-\tau(h\cdot)\Vert\leq 
\frac{\varepsilon}{2\sqrt{N}\cdot\max\Vert x_j\Vert}.
$$
Furthermore, one can
find a
finite set $F\subset S$ such that 
$$
|\tau(ux_j^*h)|\leq\frac{\varepsilon}{2\sqrt N}\quad\forall u\in
S\setminus F\quad\textrm{and}\quad\forall j=1,\ldots,N.
$$
This implies that 
$$
|\tau'(ux_j^*)|\leq\frac{\varepsilon}{\sqrt N}
\quad\forall u\in
S\setminus F\quad\textrm{and}\quad\forall j=1,\ldots,N.
$$
Let $x\in \mathrm{span}\{x_1,\ldots,x_N\}$, $\Vert x\Vert_{2,\tau'}\leq
1$. Let us write $x=\sum_{j=1}^N\xi_jx_j$, where $\xi_j=\tau'(xx_j^*)$,
and $\sum_{j=1}^N|\xi_j|^2=\Vert x\Vert_{2,\tau'}^2\leq 1$ since the $x_j$'s
are $\tau'$-orthonormal. Hence we get, for $u\in S\setminus F$:
$$
|\tau'(ux^*)| = |\sum_{j=1}^N\overline{\xi_j}\tau'(ux_j^*)|
\leq
\left(\sum_{j=1}^N|\xi_j|^2\right)^{1/2}\left(\sum_{j=1}^N|\tau'(ux_j^*)|^2\right)^{1/2}
\leq \varepsilon
$$
uniformly on the set $\{x\in\mathrm{span}\{x_1,\ldots,x_N\}\ ; \Vert
x\Vert_{2,\tau'}\leq 1\}$.\\
$(3)\ \Rightarrow\ (1)$: Let $\varphi\in M_\star$ and $\varepsilon>0$. We
 choose $x\in M$ such that 
$\Vert \varphi-\tau(\cdot x)\Vert\leq \varepsilon/2$. Applying
condition (3) to the singleton set $\{x/\Vert x\Vert_2\}$ as orthonormal set, we find a
finite subset $F$ of $S$ such that $|\tau(ux)|\leq \varepsilon/2$ for
every $u\in S\setminus F$. Hence we get $|\varphi(u)|\leq \varepsilon$
for all $u\in S\setminus F$.
\par
The last statement follows readily from condition (1).
\hfill $\square$
 
\begin{defn}
Let $M$ be a finite von Neumann algebra gifted with some fixed finite trace
$\tau$. We say that an infinite subset $S\subset U(M)$ is
\textbf{almost orthonormal} if it satisfies the equivalent conditions
in Proposition 3.2.
\end{defn}
\textbf{Remarks.} (1) Since $M$ has separable predual, an almost orthonormal
subset $S$ of $U(M)$ is necessarily countable. Indeed, let $\{x_n\
;\ n\geq 1\}$ be a $\Vert\cdot\Vert_2$-dense countable subset of the unit ball
of $M$ with respect to the operator norm. For $n\geq 1$, 
put 
$$
S_n=\{u\in S\ ;\ \max_{1\leq j\leq n}|\tau(ux_j^*)|\geq \frac{1}{n}\}.
$$
Then each $S_n$ is finite, $S_n\subset S_{n+1}$ for every $n$ and
$S=\bigcup_n S_n$.
Thus, if $S=(u_n)_{n\geq 1}$ is a sequence of unitary elements, then
$S$ is almost orthonormal if and only if $u_n$ tends weakly to
zero. In particular, every diffuse von Neumann algebra contains almost
orthonormal sequences of unitaries.
\par\noindent
(2) The reason why we choose the above definition comes from the fact that 
if $S$ is almost orthonormal in $M$, then for every $u\in S$, for
every $\varepsilon>0$, there exists a finite set $F\subset S$ such
that $|\tau(v^*u)|<\varepsilon$ for all $v\in S\setminus F$. 
A typical example of an almost orthonormal subset
in a finite von Neumann algebra is a $\tau$-orthonormal subset $S$ of 
$U(M)$ for \emph{some} trace $\tau$ on $M$: 
indeed, for every $x\in M$, the series $\sum_{u\in S}|\tau(xu^*)|^2$
converges. 
For instance, let  $v\in U(M)$ be such that $\tau(v^k)=0$
 and
for all
integers $k\in\mathbb{Z}\setminus\{0\}$. Then the subgroup generated by
$v$ is almost orthonormal.
As another example, let $\Gamma$ be a countable group and
let $\Gamma_1$ be an infinite subgroup of $\Gamma$. 
Set
$G=\lambda(\Gamma_1)$. Then $G$ is almost orthonormal in
$M=L(\Gamma)$. Indeed, if $x\in L(\Gamma)$, then
$\tau(\lambda(g)x)=\tau(x\lambda(g))=x(g^{-1})$ obviously tends to 0
as $g$ tends to infinity of $\Gamma_1$.
\par\vspace{3mm}
We come now to the main definition of our article. To motivate it,
recall that if $\Gamma$ is a (countable) group and if $\alpha$ is a
$\tau$-preserving action of $\Gamma$ on $M$, then it is called
\emph{strongly mixing} if, for every finite set $F\subset M$ and for
every $\varepsilon>0$, there exists a finite set $E\subset\Gamma$ such that
$$
|\tau(\alpha_g(a)b)-\tau(a)\tau(b)|<\varepsilon
$$
for all $a,b\in F$ and all $g\notin E$.

\begin{defn}
Let $M$ and $\tau$ be as above and let $A$ be an abelian, unital von
Neumann subalgebra of $M$. Let $G$ be a subgroup of $U(A)$.
We
say that the action of $G$ is \textbf{strongly mixing 
relative to} $A$ if, for all $x,y\in
M$, one has:
$$
\lim_{u\to\infty, u\in G}
\Vert E_A(uxu^{-1}y)-E_A(x)E_A(y)\Vert_2=0.
$$
Furthermore, we say that $A$ itself is \textbf{strongly mixing} in $M$
if, for every almost orthonormal infinite subgroup $G$ of $U(A)$, the
action of $G$ by inner automorphisms on $M$ is strongly mixing relative to $A$.
\end{defn}
\textbf{Remark.} The above property is independent
of the trace $\tau$ and it is a conjugacy invariant. Indeed, $E_A$ is
the unique conditional expectation onto $A$ and almost orthonormality
is independent of the chosen trace. 

\par\vspace{3mm}
We present now our main result.
\begin{thm}
Let $\Gamma$ be an infinite group and let $\Gamma_0$ be an infinite 
abelian
subgroup of $\Gamma$.
Let $M=L(\Gamma)$ and $A=L(\Gamma_0)$ be as above. 
Then the following conditions are equivalent:
\begin{enumerate}
\item [(1)] the action of $\Gamma_0$ by inner automorphisms on
  $M$ is strongly mixing relative to $A$;
\item [(2)] for every finite subset $C\subset\Gamma\setminus\Gamma_0$
  there exists a finite subset $E\subset\Gamma_0$ such that
  $gg_0h\notin\Gamma_0$ for all $g_0\in\Gamma_0\setminus E$ and all $g,h\in
  C$; namely, the pair $(\Gamma,\Gamma_0)$ satisfies condition (ST);
\item [(3)] for every almost orthonormal infinite subset $S\subset U(A)$,
  for all $x,y\in M$ and for every $\varepsilon>0$, there exists a
  finite subset $F\subset S$ such that
$$
\Vert E_A(uxu^*y)-E_A(x)E_A(y)\Vert_2<\varepsilon\quad \forall u\in
S\setminus F;
$$
\item [(4)] $A$ is strongly mixing in $M$.
\end{enumerate}
\end{thm}
\emph{Proof.} Trivially, $(3)\ \Rightarrow (4)\ \Rightarrow\ (1)$.\\
$(1)\ \Rightarrow\ (2)$: If $C$ is as in (2), set
$$
x=\sum_{g\in C}\lambda(g).
$$
Thus $E_A(x)=0$, and there exists a finite subset $E\subset
\Gamma_0$ such that 
$$
\Vert E_A(\lambda(g_0^{-1})x\lambda(g_0)x)\Vert_2=
\Vert E_A(x\lambda(g_0)x)\Vert_2<1\quad\forall g_0\in\Gamma_0\setminus E.
$$
But 
$$
  E_A(x\lambda(g_0)x)  = 
E_A\left(\sum_{g,h\in C}\lambda(gg_0h)\right)
  =  \sum_{g,h\in C,gg_0h\in\Gamma_0}\lambda(gg_0h),
$$
hence $\Vert E_A(x\lambda(g_0)x)\Vert_2^2=
|\{(g,h)\in C\times C\ ;\ gg_0h\in\Gamma_0\}|<1$ for all $g_0\notin
E$, which implies that $gg_0h\notin\Gamma_0$ for $g_0\notin E$ and for
all $g,h\in C$.\\
$(2)\ \Rightarrow\ (3)$: Let $S$ be an almost orthonormal infinite
subset of $U(A)$, let $x,y\in M$ and fix
$\varepsilon>0$. Using $A$-bilinearity of $E_A$, we can assume
that $E_A(x)=E_A(y)=0$, and we have to prove that one can find a
finite set $F\subset S$ such that $\Vert
E_A(uxu^*y)\Vert_2<\varepsilon$ for all $u\in S\setminus F$. To begin
with, let us assume furthermore that $x$ and $y$ have finite support,
and let us write 
$x=\sum_{g\in C}x(g)\lambda(g)$ and $y=\sum_{h\in
  C}y(h)\lambda(h)$ with $C\subset\Gamma\setminus\Gamma_0$ finite.
Let $E\subset\Gamma_0$ be as in (2) with respect to the finite set
$C$ of $\Gamma\setminus\Gamma_0$, namely
$gg_0h\notin\Gamma_0$ for all $g,h\in C$, and
$g_0\notin E$.
We claim then that $E_A(x\lambda(g_0^{-1})y)=0$
if $g_0\in\Gamma_0\setminus E^{-1}$. Indeed, if
$g_0\in\Gamma_0\setminus E^{-1}$, we have:
\begin{eqnarray*}
 E_A(x\lambda(g_0^{-1})y) & = & 
E_A\left(\sum_{g,h\in C}x(g)y(h)\lambda(gg_0^{-1}h)\right)\\
& = &  \sum_{g,h\in C,gg_0^{-1}h\in\Gamma_0}x(g)y(h)\lambda(gg_0^{-1}h)=0
\end{eqnarray*}
because $g_0^{-1}\in \Gamma_0\setminus E$.\\
Choosing $\lambda(E)\subset A$ as a $\tau$-orthonormal system, there
exists a finite subset $F$ of $S$ such that, if $u\in S\setminus F$:
$$
\sup\{|\tau(uz^*)|\ ;\ z\in \mathrm{span}(\lambda(E)),\ \Vert z\Vert_2\leq 1\}
<\frac{\varepsilon^2}{|E|\Vert x\Vert\Vert y\Vert}.
$$
Thus, for fixed $u\in S\setminus F$, take $z:=\sum_{g_0\in
  E}u(g_0)\lambda(g_0)$, so that $z\in \mathrm{span}(\lambda(E))$, $\Vert
z\Vert_2\leq 1$ and 
$$
\sum_{g_0\in E}|u(g_0)|^2=\tau(uz^*)<
\frac{\varepsilon^2}{|E|\Vert x\Vert\Vert y\Vert}.
$$
Then 
\begin{eqnarray*}
  \Vert E_A(uxu^*y)\Vert_2 & = & 
\Vert uE_A(xu^*y)\Vert_2 = \Vert E_A(xu^*y)\Vert_2\\
& \leq & 
\sum_{g_0\in E}|u(g_0)|\Vert E_A(x\lambda(g_0^{-1})y)\Vert_2<\varepsilon,
\end{eqnarray*}
using Cauchy-Schwarz Inequality.\\
Finally, if $x,y\in M$ are such that $E_A(x)=E_A(y)=0$, if
$\varepsilon>0$, let $x',y'\in L_f(\Gamma)$ be such that
$\Vert x'\Vert\leq \Vert x\Vert$, $\Vert y'\Vert\leq \Vert y\Vert$,
$E_A(x')=E_A(y')=0$ and 
$$
\Vert x'-x\Vert_2,\Vert
y'-y\Vert_2<\frac{\varepsilon}{3\cdot\max(\Vert x\Vert,\Vert y\Vert)}.
$$
Take a finite subset $F\subset S$ such that $\Vert
E_A(ux'u^*y')\Vert_2<\varepsilon/3$ for every $u\in S\setminus
F$. Then, if $u\in S\setminus F$, 
$$
\Vert E_A(uxu^*y)\Vert_2 <
\Vert y\Vert \Vert x-x'\Vert_2+\Vert x\Vert \Vert y-y'\Vert_2 + 
\frac{\varepsilon}{3}<\varepsilon.
$$
This ends the proof of Theorem 3.5.
\hfill $\square$
\par\vspace{3mm}
As in the case of weak mixing, one has for crossed products:

\begin{prop}
Let $\Gamma_0$ be an abelian group which acts
on a finite von Neumann algebra $B$ and which preserves a trace
$\tau$, then the abelian von Neumann subalgebra $A=L(\Gamma_0)$ of the
crossed product $M=B\rtimes\Gamma_0$ is strongly mixing in $M$ if and
only if the
action of $\Gamma_0$ is.
\end{prop}
In fact, Proposition 3.6 and Theorem 4.2 of \cite{Sch} prove that
weakly 
mixing MASA's are not strongly mixing in general: 
\par\vspace{3mm}\noindent
Let $\Gamma_0$ be an infinite abelian group and let $\alpha$ be a measure-preserving,
free, weakly mixing but not strongly mixing action on some standard
probability space $(X,\mathcal{B},\mu)$ as in Theorem 4.2 of \cite{Sch}. Set
$B=L^{\infty}(X,\mathcal{B},\mu)$ and let $M$ be the corresponding
crossed product $\textrm{II}_1$-factor. Then the abelian subalgebra
$A=L(\Gamma_0)$ is a weakly mixing MASA in $M$, but it is not strongly mixing.
\par\vspace{3mm}
Typical examples of strongly mixing actions are given by (generalized) 
Bernoulli shift actions: Consider a finite von Neumann algebra
$B\not=\mathbb C$ 
gifted with
some trace $\tau_B$, let
$\Gamma_0$ be an infinite abelian group that acts \emph{properly} on a countable set
$X$ :
for every finite set $Y\subset X$, the set
$\{g\in\Gamma_0\ ;\ g(Y)\cap Y\not=\emptyset\}$ is finite. Let
$(N,\tau)=\bigotimes_{x\in X}(B,\tau_B)$ be the associated infinite
tensor product. Then the corresponding Bernoulli shift action is the
action $\sigma$ of $\Gamma_0$ on $N$ given by
$$
\sigma_g(\otimes_{x\in X}b_x)=\otimes_{x\in X}b_{gx}
$$
for every $\otimes_{x}b_x\in N$ such that $b_x=1$ for all but
finitely many $x$'s. Then it is easy to see that properness of the
action implies that $\sigma$ is a strongly mixing action. The classical case
corresponds to the simply transitive action by left translations on
$\Gamma_0$.
\par\vspace{3mm}
Let $A$ be a diffuse von Neumann subalgebra of $M$ and let $Q$ be a
finite von Neumann algebra. It is proved in Theorem 2.3 of \cite{DSS}
that the normalisers $\mathcal{N}_M(A)$ and $\mathcal{N}_{M\star
  Q}(A)$ are equal. In particular, if $A$ is a singular MASA in $M$
then it is also a singular MASA in the free product $M\star Q$. 
Then, 
using the same arguments as in the proofs of Lemma 2.2, Theorem 2.3
and  Corollary
2.4 of the above mentioned article, we obtain:
\begin{prop}
Let $M$ and $Q$ be finite von Neumann algebras and let $A$ be a strongly
mixing MASA in $M$. Then $A$ is also strongly mixing in the free
product von Neumann algebra $M\star Q$. 
\end{prop}

\section{Examples}

From now on, we consider pairs $(\Gamma,\Gamma_0)$ where
$\Gamma_0$ is an abelian subgroup of $\Gamma$. We will give examples
of families of pairs that satisfy condition (ST) on the one hand, and
of pairs that satisfy condition (SS) but not (ST) on the other hand. 

\subsection{Some free products examples}

It was noted in \cite{JOF} that if $\Gamma$ is the free group $F_N$ of
rank $N\geq 2$ on free generators $a_1,\ldots, a_N$ and if $\Gamma_0$
is the subgroup generated by some fixed $a_i$, then the pair
$(F_N,\Gamma_0)$ satisfies condition (ST). See also Corollary 3.4 of \cite{SS}.
The next result extends the above case to some amalgamated products.

\begin{prop}
Let $\Gamma=\Gamma_0\star_Z \Gamma_1$ 
be an amalgamated product where $\Gamma_0$ is an infinite abelian
group, $Z$ is a finite subgroup of $\Gamma_0$ and of $\Gamma_1$.
If $\Gamma_1\not=Z$
then the pair $(\Gamma,\Gamma_0)$ satisfies
condition (ST).
\end{prop}
\emph{Proof.} Let $R$ and $S$ be sets of representatives for the left
cosets of $Z$ in $\Gamma_0$ and $\Gamma_1$ respectively, such that
$1\in R$ and $1\in S$. Recall that every element $g\in \Gamma$ has a unique
normal form
$$
g=r_1s_1\ldots r_ls_lz
$$
with $r_i\in R$, $s_i\in S$ for every $i$, 
such that
only $r_1$ or $s_1$ can be equal to 1,
and $z\in Z$. Notice that
$g$ does not belong to $\Gamma_0$ if and only if 
$s_1\not=1$.
Let $g,h\in
\Gamma\setminus\Gamma_0$. It suffices to find a finite subset
$E\subset \Gamma_0$ such that, for every $\gamma\in\Gamma_0\setminus
E$, $g\gamma h\notin\Gamma_0$. Let us write $g=r_1s_1\ldots r_ls_lz$
with $r_i\in R$, $s_i\in S$ and $z\in Z$, and $h=u_1v_1\ldots u_kv_kw$
with $u_j\in R$, $v_j\in S$ and $w\in Z$. If $\gamma=rt$ with $r\in R$
and $t\in Z$, then
\begin{eqnarray*}
g\gamma h &=& r_1s_1\ldots r_ls_lzrtu_1v_1\ldots u_kv_kw  \\
 & = & r_1s_1\ldots r_ls_l(ru_1)v'_1z'_1\ldots u_kv_kw
\end{eqnarray*}
with $v'_1\in S$ and $z'_1\in Z$ such that $v'_1z'_1=(zt)v_1$ is the
decomposition of $(zt)v_1$ in the partition $\Gamma_1=\coprod_{\sigma\in
  S}\sigma Z$.
Observe that $v'_1\not=1$ because $v_1\not=1$. Continuing in the same
way, we move elements of $Z$ to the right as most as possible and we get finally
$$
g\gamma h=r_1s_1\ldots r_ls_l(ru_1)v'_1\ldots u_kv'_kw'
$$
with every $v'_j\in S$ and $w'\in Z$.
If $s_l\not=1$
one can take $E=Z\cup u^{-1}_1Z$, and if $s_l=1$ one can take
$E=Z\cup (r_lu_1)^{-1}Z$.
\hfill $\square$
\par\vspace{3mm}
Finally, we get from Theorem 3.5 and Proposition 3.7:

\begin{prop}
Let $\Gamma$ be an infinite group, let $\Gamma_0$ be an infinite
abelian subgroup of $\Gamma$ such that the pair $(\Gamma,\Gamma_0)$
satisfies condition (ST) and let $G$ be any countable group. Then the
pair $(\Gamma\star G,\Gamma_0)$ satisfies also condition (ST).
\end{prop}


\subsection{The case of malnormal subgroups}

In \cite{Popa83},
S. Popa introduced a condition on pairs $(\Gamma,\Gamma_0)$
 in order to obtain orthogonal
pairs of von Neumann subalgebras; it was also used later in 
\cite{SS} to get asymptotic homomorphism
conditional expectations. We say that $\Gamma_0$ is a
\textbf{malnormal} subgroup of $\Gamma$ if it satisfies the following condition:
\par
\emph{$(\star)$ For every $g\in \Gamma\setminus \Gamma_0$, 
one has $g\Gamma_0g^{-1}\cap\Gamma_0=\{1\}$.}
\par\vspace{3mm}
Then we observe that condition $(\star)$ implies condition (ST). Indeed, for $g,h\in
\Gamma$, set $E(g,h)=\{\gamma\in \Gamma_0\ ;\ g\gamma h\in\Gamma_0\}=
g^{-1}\Gamma_0h^{-1}\cap\Gamma_0$. If $(\Gamma,\Gamma_0)$
satisfies $(\star)$, and if $g,h\in\Gamma\setminus\Gamma_0$, then $E(g,h)$
contains at most one element (see the proof of Lemma 3.1 of
\cite{SS}).
In turn, condition (ST) means exactly that $E(g,h)$ is finite for all 
$g,h\in\Gamma\setminus\Gamma_0$. Thus,
if $\Gamma_0$ is torsion free, then conditions $(\star)$ and (ST) are
equivalent because, in this case, if $g\in\Gamma\setminus\Gamma_0$,
the finite set $E(g^{-1},g)=g\Gamma_0g^{-1}\cap\Gamma_0$ is a finite
subgroup of $\Gamma_0$, hence is trivial. 
\par\vspace{3mm}
However, condition $(\star)$ is strictly stronger than condition (ST)
in general: let $\Gamma=\Gamma_0\star_Z \Gamma_1$ be as in Proposition
4.1 above, and assume further that $Z\not=\{1\}$ and that there exists
$g\in\Gamma_1\setminus Z$ such that $zg=gz$ for every $z\in Z$. Then 
$g\Gamma_0 g^{-1}\cap\Gamma_0\supset Z$, and $(\Gamma,\Gamma_0)$ does
not satisfy $(\star)$. Observe that, in this case, $\Gamma$ is not
necessarily an ICC group. 
However, replacing it by 
a non trivial free product group
$\Gamma\star G$, we get a pair $(\Gamma\star G,\Gamma_0)$ satisfying
condition (ST) by Proposition 3.7
but not $(\star)$, and $\Gamma\star G$ is an ICC group.
\par\vspace{3mm}
It is known that in some classes of groups, maximal abelian subgroups
are malnormal. 
This is e.g. the case in hyperbolic groups \cite{Gro} or in groups acting freely on $\Lambda$-trees \cite{Bas}.
We present here explicitely a sufficient condition to get malnormal subgroups in groups acting on trees:

\begin{prop}\label{tree}
 Let $\Gamma$ be a group acting on a tree $T$ without inversion. 
Let $t\in \Gamma$ and set $\Gamma_0 = \langle t \rangle$. 
Assume that:
 \begin{enumerate}
 \item [(1)] there exist neither $u\in\Gamma$ nor $n\in\mathbb
   Z\setminus\{\pm 1\}$ such
   that $t=u^n$;
 \item [(2)] the induced automorphism of $T$ (again denoted by $t$) is hyperbolic;
 \item [(3)] the subgroup $\bigcap_{v\in \axis(t)} \Stab(v)$ is trivial;
 \item [(4)] no element of $\Gamma$ induces a reflection of $\axis(t)$.
 \end{enumerate}
Then $\Gamma_0$ is a malnormal subgroup of $\Gamma$, and in particular
the pair $(\Gamma,\Gamma_0)$ satisfies condition (ST).
\end{prop}

Before proving this, we introduce the necessary terminology about
actions on trees. 
If a group $G$ acts without inversion on a tree $X$, it is well-known
\cite[Chap I.6.4]{Ser} 
that an element $g\in G$ is either \emph{elliptic}, that is it fixes a
vertex, or \emph{hyperbolic}, 
that is $g$ preserves an infinite geodesic (called \emph{axis}) on
which it acts by a non-trivial translation. 
It is easy to see that a hyperbolic element $g$ has a unique axis. It will be denoted $\axis(g)$.

\par\vspace{3mm}
\emph{Proof of Proposition \ref{tree}.}  Given an element $g\in \Gamma$ such that $g\Gamma_0
g^{-1} \cap \Gamma_0 \neq\{1\}$, we have to prove that $g\in
\Gamma_0$. 
First, we write $gt^kg^{-1} = t^{k'}$ for some $k,k' \in \mathbb{Z}^*$. The
elements $t^k$ and $t^{k'}$ having $\textrm{axis}(t)$ as axis, $g$ has to
preserve $\axis(t)$. 
Then, by hypotheses (3) and (4), $g$ is an hyperbolic element and
$\axis(g) = 
\axis(t)$.
 Let us denote by $\ell(\gamma)$ the translation length of an element
 $\gamma\in\Gamma$. 
By B\'ezout's Theorem, there exists an element $s=g^m t^n$ with
 $m,n\in\mathbb Z$ such that $\ell(s)$ is the greatest common divisor of
 $\ell(g)$ and $\ell(t)$. 
Then, there exist $\delta,\varepsilon \in\{\pm 1\}$ such that
 $gs^{\delta \ell(g)/\ell(s)}$ 
and $ts^{\varepsilon \ell(t)/\ell(s)}$ fix every vertex of
 $\axis(t)$. By hypothesis (3), we get $g=s^{-\delta \ell(g)/\ell(s)}$ 
and $t=s^{-\varepsilon \ell(t)/\ell(s)}$. Then (1) gives $t=s^{\pm
 1}$, so that $g=t^{\pm \ell(g)/\ell(t)}$. 
As desired, $g$ is an element of $\Gamma_0$.
\hfill $\square$
\par\vspace{3mm}
As it will be seen below, wide families of HNN extentions satisfy
hypotheses of Proposition 4.3. Thus, 
let $\Gamma=HNN(\Lambda,H,K,\phi)$ be an HNN extension where $H,K$ are
 subgroups of $\Lambda$ and, as usual, where $\phi:H\rightarrow
K$ is an isomorphism. Denote by $t$ the stable letter such that
$t^{-1}ht=\phi(h)$ for all $h\in H$, and by $\Gamma_0$ the subgroup 
generated by $t$. Recall that a sequence
$g_0,t^{\varepsilon_1},\ldots,t^{\varepsilon_n},g_n$, $(n\geq 0)$ is
\emph{reduced} if $g_i\in \Lambda$ and 
$\varepsilon_i=\pm 1$ for every $i$, and if there is no subsequence
$t^{-1},g_i,t$ with $g_i\in H$ or $t,g_i,t^{-1}$ with $g_i\in K$. As
is well known, if the sequence 
$g_0,t^{\varepsilon_1},\ldots,t^{\varepsilon_n},g_n$ is reduced and if
$n\geq 1$ then the corresponding element 
$g=g_0t^{\varepsilon_1}\cdots t^{\varepsilon_n}g_n\in\Gamma$ is non trivial
 (Britton's
Lemma).
We also say that such an element is in \emph{reduced form}.
Furthermore, if 
$g=g_0t^{\varepsilon_1}\cdots t^{\varepsilon_n}g_n$ and
$h=h_0t^{\delta_1}\cdots t^{\delta_m}g_m$ are in reduced form and if $g=h$, then
$n=m$ and $\varepsilon_i=\delta_i$ for every $i$. Hence the
\emph{length} $\ell(g)$ of 
$g=g_0t^{\varepsilon_1}\cdots t^{\varepsilon_n}g_n$ (in reduced form) 
is the integer $n$. Finally, recall from \cite{YS} that, for
every positive integer $j$, 
$\mathrm{Dom}(\phi^j)$ is defined by $\mathrm{Dom}(\phi)=H$
for $j=1$ and, by induction, 
$\mathrm{Dom}(\phi^j)=\phi^{-1}(\mathrm{Dom}(\phi^{j-1})\cap K)\subset
H$ for $j\geq 2$.
\par\vspace{3mm}
If $\Gamma=HNN(\Lambda,H,K,\phi)$, its \emph{Bass-Serre tree} has
$\Gamma/\Lambda$  
as set of vertices and $\Gamma/H$ as set of oriented edges. The origin
of the edge 
$\gamma H$ is $\gamma \Lambda$ and its terminal vertex is $\gamma t
\Lambda$. 
Chapter I.5 in \cite{Ser}, and Theorem 12 in particular, ensures that
it is a tree. 
It is obvious that the Bass-Serre tree is endowed with an orientation-preserving $\Gamma$-action.

\begin{cor}\label{HNNtree}
Suppose that for each $\lambda \in \Lambda^*$, there exists $j>0$ such
that $\lambda 
\notin \Dom(\phi^j)$.
Then $\Gamma_0$ is a malnormal subgroup of $\Gamma$.
\end{cor}
\emph{Proof.}  Let $T$ be the Bass-Serre tree of the HNN-extension. We
check the hypotheses of Proposition 4.3. 
Since (1), (2) and (4) are obvious, we prove (3).
  We have $\bigcap_{v\in \axis(t)} \Stab(v) =
  \bigcap_{k\in\mathbb Z} \Lambda_k$, where $\Lambda_k = t^{-k}\Lambda t^k$. 
Assume by contradiction that the intersection contains a non trivial
element $\lambda$. By hypothesis, there exists $j\in\mathbb{N}^*$, 
such that $\lambda\not\in \Dom(\phi^j)$ and we may
assume $j$ to be minimal for this property. 
This means that $t^{1-j}\lambda t^{j-1} = \phi^{j-1}(\lambda) \in
\Lambda\setminus H$, and $t^{-j}\lambda t^{j} \not\in \Lambda$ 
by Britton's Lemma. We get $\lambda \not \in \Lambda_{-j}$, a
contradiction. 
This proves (3).
\hfill $\square$
\par\vspace{3mm}
As it will be recalled in the first example below,
Thompson's group $F$ is an HNN extention, and it satisfies condition
(ST) with respect to the subgroup generated by $x_0$, by Lemma 3.2 of \cite{JOF}.
\par\vspace{3mm}\noindent
\textbf{Examples.} (1) 
For every integer $k\geq 1$, denote by $F_k$ the subgroup of $F$
generated by $x_k,x_{k+1},\ldots$, and denote by $\sigma$ the ``shift
map'' defined by $\sigma(x_n)=x_{n+1}$, for $n\geq 0$. Its restriction
to $F_k$ is an isomorphism onto $F_{k+1}$, and in particular, the
inverse map $\phi:F_2\rightarrow F_1$ is an isomorphism which
satisfies $\phi(x)=x_0xx_0^{-1}$ for every $x\in F_2$. As in
Proposition 1.7 of \cite{BG}, it is evident that $F$ is the 
HNN extension $HNN(F_1,F_2,F_1,\phi)$
with $t=x_0^{-1}$ as stable
letter. With these choices, $F$ satisfies the hypotheses of
Proposition 4.1, and this proves that the pair $(F,\Gamma_0)$
satisfies condition (ST).\\
(2) Let $m$ and $n$ be non-zero integers. The associated
\emph{Baumslag-Solitar group} is the group which has the following
presentation:
$$
BS(m,n)=\langle a,b\ |\ ab^ma^{-1}b^{-n}\rangle.
$$
Since $a^{-1}b^na=b^m$, $BS(m,n)$ is an HNN extension
$HNN(\mathbb{Z},n\mathbb{Z},m\mathbb{Z},\phi)$ where $\phi(nk)=mk$ for
every integer $k$. Assume first that $|n|>|m|$ and denote by
$\Gamma_0$ the subgroup generated by $a$.
Then it is easy to
check that the pair $(\Gamma_0,BS(m,n))$ satisfies the condition in
Proposition 4.1. Thus, it satisfies also condition (ST). If $|m|>|n|$,
replacing $a$ by $a^{-1}$ (which does not change $\Gamma_0$), one gets
the same conclusion. Thus, all Baumslag-Solitar groups $BS(m,n)$ with
$|m|\not=|n|$ satisfy condition (ST) with respect to the subgroup
generated by $a$. Observe that the latter class is precisely the class
of Baumslag-Solitar groups that are ICC (\cite{YS}).
\par\vspace{3mm}
We turn to free products. The \emph{Bass-Serre tree} of $\Gamma=A*B$
has 
$(\Gamma/A) \sqcup (\Gamma/B)$ as set of vertices and $\Gamma$ as set
of oriented edges. The origin of the 
edge $\gamma$ is $\gamma A$ and its terminal vertex is $\gamma
B$. Again, Chapter I.5 in \cite{Ser} 
ensures that it is a tree. Here is another consequence of Proposition 4.3:

\begin{cor}
 Let $\Gamma = A*B$ be a free product such that $|A|\geq 2$ and $B$
 contains an element $b$ of order at least $3$. 
Let $a$ be a non trivial element of $A$ and $\Gamma_0 = \langle ab
 \rangle$. 
Then $\Gamma_0$ is malnormal in $\Gamma$.
\end{cor}
\emph{Proof.}
 Let $T$ be the Bass-Serre tree of the free product. Again, we check
 the hypotheses of Proposition 4.3, with $t=ab$. 
Since (1) and (2) are obvious, we focus on (3) and (4).
 The axis of $t$ has the following structure.
$$
\ldots \leftarrow b^{-1}a^{-1}b^{-1}A \rightarrow  b^{-1}a^{-1}B  \leftarrow
 b^{-1}A \rightarrow  B \leftarrow A \rightarrow  aB \leftarrow
 abA\rightarrow \ldots
$$
 
 To prove (3), we remark that the vertices $A$ and $abA$ are on
 $\axis(t)$ and that $\operatorname{Stab}(A) \cap
 \Stab(abA) = 
A \cap abAb^{-1}a^{-1} = \{1\}$.\\
 To prove (4), assume by contradiction that there exists an elliptic
 element $g\in\Gamma$ which induces a reflection on $\axis(t)$. 
Up to conjugating $g$ by a power of $t$, we may assume that the fixed 
point of $g$ on $\axis(t)$ is $A$ or $B$. 
 If $g$ fixes $A$, we have $gb^{-1}A = abA$, so that $g\in A$ and
 $bg^{-1}ab \in A$. This is a contradiction since $|b|\geq 3$.
 Now, if $g$ fixes $B$, we have $gabA = b^{-1}a^{-1}b^{-1}A$, so that
 $g\in B$ and $babgab \in A$. This implies $g=b^{-1}$, $a^2=1$, and
 $b^2=1$, 
which is again a contradiction since $|b|\geq 3$.
\hfill $\square$


\subsection{Some semidirect products}

Let $H$ be a discrete group, let $\Gamma_0$ be an infinite
abelian group and let
$\alpha:\Gamma_0\rightarrow\mathrm{Aut}(H)$ be an action of $\Gamma_0$ on $H$. Then
the semi-direct product $\Gamma=H\rtimes_{\alpha}\Gamma_0$ is the direct
product set $H\times\Gamma_0$ gifted with the multiplication
$$
(h,\gamma)(h',\gamma')=(h\alpha_\gamma(h'),\gamma\gamma')\quad
\forall (h,\gamma),(h',\gamma')\in \Gamma.
$$
The action $\alpha$ lifts from $H$ to the von Neumann algebra $L(H)$,
and $L(\Gamma)=L(H)\rtimes \Gamma_0$ is a crossed product von Neumann
algebra in a natural way. In Theorem 2.2 of \cite{RSS}, the authors consider a
sufficient condition on the action $\alpha$ on $H$ which ensures that
$L(\Gamma_0)\subset L(\Gamma)$ is a strongly singular MASA in $L(\Gamma)$ and
that $L(\Gamma)$ is a type $\textrm{II}_1$ factor. In
fact, it turns out that their condition implies that $L(\Gamma_0)$ is
strongly mixing in $L(\Gamma)$, as we prove here:
\begin{prop}
Let $H$ and $\Gamma_0$ be infinite discrete groups, $\Gamma_0$ being abelian, 
let $\alpha$ be an action of $\Gamma_0$ on $H$ and let 
$\Gamma=H\rtimes_{\alpha}\Gamma_0$. If, for each $\gamma\in \Gamma_0^*$, the only fixed
point of $\alpha_\gamma$ is $1_H$, then:
\begin{enumerate}
\item [(1)]  $\Gamma$ is an ICC group;
\item [(2)]  the
pair $(\Gamma,\Gamma_0)$ satisfies condition (ST);
\item [(3)] the action of $\Gamma_0$ on $L(H)$ is strongly mixing.
\end{enumerate}
 In particular, $L(\Gamma)$
is a type $II_1$ factor and
$L(\Gamma_0)$ is
strongly mixing in $L(\Gamma)$.
\end{prop}
\emph{Proof.}  Statement (1)
is proved in \cite{RSS}. 

Thus it remains to prove (2) and (3).
\par
To do that, we claim first that if $\alpha$ is as 
above, then the triple $(\Gamma_0,H,\alpha)$ satisfies the following
condition whose proof is inspired by that of Theorem 2.2 of \cite{RSS}:
\par\vspace{3mm}
\emph{For every finite subset $F$ of $H^*$, there exists a finite set
  $E$ in $\Gamma_0$ such that $\alpha_\gamma(F)\cap F=\emptyset$ for all $\gamma\in
  \Gamma_0\setminus E$.}
\par\vspace{3mm}
Indeed, if $F$ is fixed, set $I(F)=\{\gamma\in \Gamma_0\ ;\ 
\alpha_\gamma(F)\cap F\not=\emptyset\}$.
If $I(F)$ would be infinite for some $F$, set for $f\in F$:
$$
S_f=\{\gamma\in \Gamma_0\ ;\ \alpha_\gamma(f)\in F\},
$$
so that $I(F)=\bigcup_{f\in F}S_f$, and $S_f$ would be infinite for at
least one $f\in F$. There would exist then distinct elements
$\gamma_1$ and 
$\gamma_2$ 
of $\Gamma_0$ such that $\alpha_{\gamma_1}(f)=\alpha_{\gamma_2}(f)$, since $F$ is
finite. This is impossible because $\alpha_{\gamma_1^{-1}\gamma_2}$ cannot have
any fixed point in $F\subset H^*$. 
\par
Let us prove (2).
Fix some finite set $C\subset \Gamma\setminus \Gamma_0$. Without loss of
generality, we assume that $C=C_1\times C_2$ with $C_1\subset H^*$ and
$C_2\subset \Gamma_0$ finite. Take $F=C_1\cup C_1^{-1}\subset H^*$ in the
condition above and let $E_1\subset \Gamma_0$ be a finite set such that
$\alpha_{\gamma}(F)\cap F=\emptyset$ for all $\gamma\in \Gamma_0\setminus E_1$. Put 
$E=\bigcup_{\gamma\in C_2}\gamma^{-1}E_1$, which is finite. Then it is easy to
check that for all $(h,\gamma),(h',\gamma')\in C$ and for every  
$g\in \Gamma_0\setminus E$, one has
$$
(h,\gamma)(1,g)(h',\gamma')=
(h\alpha_{\gamma g}(h'),\gamma g\gamma')\notin \{1_H\}\times \Gamma_0.
$$
Finally, in order to prove (3), it suffices to see that, if $a,b\in
L_f(H)$, then there exists a finite subset $E$ of $\Gamma_0$ such that 
$$
\tau(\alpha_\gamma(a)b)=\tau(a)\tau(b)\quad\forall \gamma\notin E.
$$
Thus, let $S\subset \Gamma_0$ be a finite subset such that 
$a=\sum_{\gamma\in S}a(\gamma)\lambda(\gamma)$ and $b=\sum_{\gamma\in
  S}b(\gamma)\lambda(\gamma)$.
Choose $E\subset\Gamma_0$ finite such that 
$\alpha_{\gamma}(S^*\cup (S^*)^{-1})\cap(S^*\cup(S^*)^{-1})=\emptyset$
for every $\gamma\in\Gamma_0\setminus E$. Then, if $\gamma\notin E$,
we have
$$
\tau(\alpha_\gamma(a)b)=\sum_{h_1,h_2\in S}a(h_1)b(h_2)
\tau(\lambda(\alpha_\gamma(h_1)h_2))=\tau(a)\tau(b)
$$
since $\alpha_\gamma(h_1)h_2\not=1$ for $h_1,h_2\in S^*$.
 \hfill $\square$ 
\par\vspace{3mm}\noindent
\textbf{Example.} Let $d\geq 2$ be an integer and let $g\in GL(d,\mathbb
Z)$. Then it
defines a natural action of $\Gamma_0=\mathbb Z$ on $H=\mathbb{Z}^d$
which has no non trivial fixed point if and only if 
the list of eigenvalues of $g$ contains no root of unity. (See for
instance Example 2.5 of \cite{YS}.) 
\par\vspace{3mm}
We give below an application of Proposition 4.6 to some HNN
extensions (a wider class than in Corollary \ref{HNNtree});
as it follows from \cite[Chap I.1.4, Prop 5]{Ser}, if $\Gamma = HNN(\Lambda,H,K,\phi)$ with
stable letter $t$ 
such that $t^{-1}ht = \phi(h)$ for all $h\in H$, if  $\Gamma_0=
\langle t \rangle$ as above and if $\sigma_t$ denotes the quotient map
$\Gamma\rightarrow \Gamma_0$, then $\Gamma$ is a semidirect product
group $N\rtimes \Gamma_0$ where $N=\operatorname{ker}(\sigma_t)$ is a
direct limit, or \emph{amalgam}, of the system
$$
\ldots \Lambda_{-1} \nwarrow K_{-1}=H_0\nearrow\Lambda_0\nwarrow
K_0=H_1 \nearrow \Lambda_1 \ldots
$$
where $\Lambda_i = t^{-i}\Lambda t^i$, $H_i = t^{-i}H t^i$ and $K_i =
t^{-i}K t^i$ for all $i\in \mathbb{Z}$. 
The conjugation operation $n \mapsto t^{-1}n t$ corresponds 
then to a shift to the right direction.

\begin{cor}
Suppose that for all $j\in\mathbb{N}^*$, the homomorphism $\phi^j$ has no non
trivial fixed point, 
that is, for all $h\in H$, $\phi^j(h) = h$ implies $h=1$.
Then the following hold:
\begin{enumerate}
 \item the group $\Gamma$ is ICC;
 \item the pair $(\Gamma,\Gamma_0)$ satisfies condition (ST);
 \item the algebra $L(\Gamma)$ is a type II$_1$ factor, in which $L(\Gamma_0)$ is strongly mixing.
\end{enumerate}
\end{cor}
\emph{Proof.} It suffices to prove that, for any $k\in\mathbb{Z}^*$
and any $n\in N^*$, one has $t^{-k}nt^k \not= n$.\\
 Assume by contradiction that $t^{-k}nt^k = n$ for some $k\in \mathbb{Z}^*$,
 and some $n \in N^*$. 
Up to replacing $k$ by $-k$, we assume that $k\in\mathbb{N}^*$. Then there exists
$s\in \mathbb{N}$ such that $n$ is in the subgroup of $N$ generated by 
$\Lambda_{-s}, \ldots, \Lambda_{0}, \ldots, \Lambda_{s}$. Let now
$\ell$ be a multiple of $k$ such that $\ell-s> s+k$. 
The element $t^{-\ell}nt^\ell$ is in the subgroup of $N$ generated by
$\Lambda_{\ell-s}, \ldots, \Lambda_{\ell}, \ldots, \Lambda_{\ell+s}$. 
Then we consider the subgroups $N_l$ generated by $\ldots,
\Lambda_{s-1}, \Lambda_{s}$, and $N_r$ generated 
by $\Lambda_{s+1}, \Lambda_{s+2}, \ldots$, and, since $\ell-s > s$, 
we have $N = N_l *_{K_s=H_{s+1}} N_r$ with $n\in N_l$ and
$t^{-\ell}nt^\ell \in N_r$. Since $\ell$ is a multiple of $k$, 
we have $t^{-\ell}nt^\ell = n$ and this element is in $K_s =
H_{s+1}$. 
By similar arguments (shifting the ``cutting index'' in the
construction of $N_l, N_r$) we obtain that $n \in K_{s+1} = H_{s+2}$,
\ldots, 
$n \in K_{\ell-s-1} = H_{\ell-s}$. Hence, the elements $n':=
t^{\ell-s} n t^{s-\ell}$, $t^{-1}n't= t^{\ell-s-1} n t^{s-\ell+1}$, 
\ldots, $t^{-\ell+2s+1}n't^{\ell-2s-1}= t^{s+1} n t^{-s-1}$ are in
$H_0=H$. Thus, since $\ell - s > k$, we have $n', t^{-1}n't, 
\ldots, t^{-k}n't^k \in H$. Consquently, $\phi^k(n')$ exists, and
$\phi^k(n') = t^{-k}n't^k = t^{\ell-s-k}nt^{k+s-\ell} = 
t^{\ell-s}n't^{s-\ell} = n'$. On the other hand, $\phi^k$ 
has no non trivial fixed point. This is a contradiction.
\hfill $\square$


\subsection{Final examples and remarks}

Next, let us look at examples where $\Gamma_0$ is not cyclic (inspired by
Sinclair and Smith, \cite{SSP}):
let $\mathbb Q$ 
be the additive group of rational numbers and denote
by $\mathbb{Q}^\times$ the multiplicative group of nonzero rational
numbers.

For each positive integer $n$, set
$$
F_n=\{\frac{p}{q}\cdot 2^{kn}\ ; \ p,q\in\mathbb{Z}_{\mathrm{odd}},
\ k\in\mathbb{Z}\}\subset \mathbb{Q}^\times
$$
and 
$$
F_\infty=\{\frac{p}{q}\ ; \ p,q\in\mathbb{Z}_{\mathrm{odd}}\}
\subset \mathbb{Q}^\times.
$$
Next, for $n\in\mathbb{N}\cup\{\infty\}$, set
$$
\Gamma(n)=\left\lbrace
\left(
\begin{array}{cc}
f & x\\
0 & 1
\end{array}\right)\ ;\ 
f\in F_n,\ x\in\mathbb{Q}
\right\rbrace
$$
and let $\Gamma_0(n)$ be the subgroup of diagonal elements of 
$\Gamma(n)$. $\Gamma(n)$ is an ICC,
amenable group. Then the pair $(\Gamma(n),\Gamma(n)_0)$
satisfies condition (ST) for every $n$. 

However, if we consider larger matrices, the corresponding pairs of groups
do not satisfy condition (ST). 
Let us fix two positive integers $m$ and $n$, and set
$$
\Gamma(m,n)=
\left\lbrace
\left(
\begin{array}{ccc}
1 & x & y \\
0 & f_1 & 0 \\
0 & 0 & f_2
\end{array}\right)\ ;\ 
f_1\in F_m,\ 
f_2\in F_n,\ x,y\in\mathbb{Q}
\right\rbrace
$$
and let $\Gamma_0(m,n)$ be the corresponding diagonal subgroup.
Then we have:
$$
\left(
\begin{array}{ccc}
1 & 0 & 1 \\
0 & 1 & 0 \\
0 & 0 & 1
\end{array}\right)
\left(
\begin{array}{ccc}
1 & 0 & 0 \\
0 & f & 0 \\
0 & 0 & 1
\end{array}\right)
\left(
\begin{array}{ccc}
1 & 0 & -1 \\
0 & 1 & 0 \\
0 & 0 & 1
\end{array}\right)=
\left(
\begin{array}{ccc}
1 & 0 & 0 \\
0 & f & 0 \\
0 & 0 & 1
\end{array}\right)
$$
which belongs to $\Gamma_0$
for all $f\in F_m$. Thus the pair $(\Gamma(m,n),\Gamma(m,n)_0)$
does not satisfy condition (ST), though it satisfies condition (SS).

\bibliographystyle{plain}
\bibliography{refstronglysing}
\par
\vspace{1cm}
\noindent
       P. J.\\
       Universit\'e de Neuch\^atel,\\
       Institut de Math\'emathiques,\\       
       Emile-Argand 11\\
       Case postale 158\\
       CH-2009 Neuch\^atel, Switzerland\\
       \small{paul.jolissaint@unine.ch}

       \par\vspace{3mm}\noindent
       Y. S.\\
       Laboratoire de Math\'ematiques,\\
       Universit\'e Blaise Pascal\\
       Campus universitaire des C\'ezeaux\\
       63177 Aubi\`ere Cedex, France\\
       \small{yves.stalder@math.univ-bpclermont.fr}

\end{document}